\documentclass[final,1p,times]{elsarticle}
\usepackage{amssymb}
\usepackage{amsfonts}
\usepackage{amsmath}
\usepackage{geometry}
\usepackage{graphicx}
\usepackage{epsfig}
\usepackage{epstopdf}
\usepackage{bm}
\numberwithin{equation}{section}
\newtheorem{thm}{Theorem}

\newdefinition{rmk}{Remark}

\newproof{pf}{Proof} 

\newproof{pot}{Proof of Theorem  
}
\begin{document}
\begin{frontmatter}
\title{Impossibility of convergence of a Heun function on the boundary of the disc of convergence}
\author{Yoon-Seok Choun}
\ead{ychoun@gradcenter.cuny.edu; ychoun@gmail.com}
\address{Department of Physics, Hanyang University, Seoul, 133-791, South Korea}
\begin{abstract}

The Heun's equation is the Fuchsian equation of second order with four regular singularities. 
Heun functions generalize well-known special functions such as Spheroidal Wave, Lam\'{e}, Mathieu, hypergeometric-type functions, etc.  
The recursive relation of coefficients starts to appear by putting a power series into the Heun equation. A local Heun function consists of the three term recurrence relation in its power series, and we prove that the function is not convergent on the boundary of the disc of convergence. 
\end{abstract}

\end{frontmatter}
\section{Introduction}\label{sec.1}

A homogeneous linear ordinary differential equation of order $j$ with variable coefficients is of the form
\begin{equation}
a_j(x)y^{(j)}+ a_{j-1}(x)y^{(j-1)}+ \cdots + a_1(x)y^{'} + a_0(x)y =0
 \label{jj:1}
\end{equation}
Assuming its solution as a power series in the form
\begin{equation}
y(x) = \sum_{n=0}^{\infty } d_n x^{n}
 \label{jj:2}
\end{equation}
where $d_0 =1$ chosen for simplicity from now on. For a fixed $k\in \mathbb{N}$,  we suggest that we arrive at the $(k+1)$-term recurrence relation putting (\ref{jj:2}) in (\ref{jj:1}).
\begin{equation}
d_{n+1} = \alpha _{1,n} d_{n} + \alpha _{2,n} d_{n-1} + \alpha _{3,n} d_{n-2} + \cdots + \alpha _{k,n} d_{n-k+1}
 \label{jj:3}
\end{equation}
with seed values $d_{j} = \sum_{i=1}^{j}\alpha _{i,j-1} d_{j-i}  $ where $j=1,2,3,\cdots, k-1$ \& $k\in \mathbb{N}-\{1\}$.
\begin{thm}
Assuming ${\displaystyle \lim_{n\rightarrow \infty}\alpha _{l,n}= \alpha _{l}<\infty }$ in (\ref{jj:3}),
the domain of absolute convergence of (\ref{jj:2}) is written by \cite{Choun2018a}
\begin{equation}
\mathcal{D} :=  \left\{  x \in \mathbb{C} \Bigg| \sum_{m=1}^{k}\left| \alpha _{m} x^m \right| <1 \right\}
 \label{jj:4}
\end{equation}
\end{thm} 
For a 2-term case,  the radius of convergence is given by putting $k=1$ in (\ref{jj:4}). 
One of well-known examples of the 2-term case is a hypergeometric function, and its convergent domain is taken by
\begin{equation}
\mathcal{D} :=  \left\{  x \in \mathbb{C} \big|  \lim_{n\rightarrow \infty}\left| \frac{(n+a)(n+b)}{(n+c)(n+1)} x\right|  = \left|  x \right| <1 \right\} \nonumber
\end{equation}
where $a,b,c \in \mathbb{C}$. 

What happens when $|x|$ equals to the radius of convergence (at $|x|=1$ for the hypergeometric series)? 
The series might converge at the boundary (at both endpoints for real $x$), diverge at both, or converge only one of these values. If it converges at both, the convergence might be absolute or conditional. In 1812 Carl Friedrich Gauss published a test that determines that convergence at the endpoints for every power series \cite{Birk,Gauss}. It is a definitive test that works when the power series is hypergeometric. For the hypergeometric function, a series is absolute convergent as $\mathbb{R}(c) > \mathbb{R}(a+b)$.

However, more than a 3-term case, it has been unknown what values coefficients make a series as absolute convergent at both endpoints. We work on the problem how  a power series behaves on the boundary of the disc of convergence for the three term recurrence relation (including a local Heun function).
 In other words, what condition makes a series as absolutely convergent at $ \sum_{m=1}^{2}\lim_{n\rightarrow \infty}\left| \alpha _{m,n} x^m\right| =1 $?   
Surprisingly, unlike the 2-term case such as a Hypergeometric function, to begin with, a local Heun function which is the 3-term case does not converge.
In this paper, we prove that  the function is not convergent on the boundary of the disc of convergence.

The Heun's equation is a Fuchsian second-order differential equation of the form \cite{Heun1889,Ronv1995}
\begin{equation}
\frac{d^2{y}}{d{x}^2} + \left(\frac{\gamma }{x} +\frac{\delta }{x-1} + \frac{\varepsilon  }{x-a}\right) \frac{d{y}}{d{x}} +  \frac{\alpha \beta x-q}{x(x-1)(x-a)} y = 0 \label{eq:1}
\end{equation}
with the condition $\varepsilon = \alpha +\beta -\gamma -\delta +1$ and $a \ne 0 $. It has four regular singular points which are 0, 1, $a$ and $\infty $ with exponents $\{ 0, 1-\gamma \}$, $\{ 0, 1-\delta \}$, $\{ 0, 1-\varepsilon \}$ and $\{ \alpha, \beta \}$. Assume that $y(x)$ has a series expansion of the form
\begin{equation}
y(x)= \sum_{n=0}^{\infty } d_n x^{n+\lambda } \label{eq:2}
\end{equation}
where $\lambda $ is an indicial root. Plug (\ref{eq:2})  into (\ref{eq:1}):
\begin{equation}
d_{n+1}=A_n \;d_n +B_n \;d_{n-1} \hspace{1cm};n\geq 1 \label{eq:3}
\end{equation}
with $A_n= A\; \overline{A}_n$ and $B_n= B\; \overline{B}_n$
\begin{equation}
\begin{cases}
A = \frac{1+a}{a} \\
B = -\frac{1}{a} \\
\overline{A}_n  = \dfrac{n^2 +\dfrac{\alpha +\beta -\delta +2\lambda +a(\gamma +\delta -1+2\lambda )}{1+a} n + \dfrac{\lambda \left( \alpha +\beta -\delta +\lambda +a (\gamma +\delta -1+\lambda )\right) +q}{1+a}}{n^2 +(1+\gamma +2\lambda )n +(1+\lambda )(\gamma +\lambda )} \\
\overline{ B}_n   = \dfrac{n^2 +(\alpha +\beta -2+2\lambda )n +(\alpha -1+\lambda ) (\beta -1+\lambda )}{n^2 + (\gamma +1+2\lambda )n+(1+\lambda )(\gamma +\lambda )} \\
d_1= A_0 \;d_0 = A\; \overline{ A_0} \;d_0
\end{cases}\label{eq:4}
\end{equation}
We have two indicial roots which are $\lambda = 0$ and $ 1-\gamma $. And a local Heun function of (\ref{eq:1}) around $x=0$ is absolutely convergent where $\left| A x\right| + | B x^2 |<1 $ \cite{Choun2013}. 
%
%
%
%
\section{Result}\label{sec.2}
From a methodological perspective for proof, we study the nature of a Heun function by utilizing Gauss's Test (which is suitable for a 3-term recurrence relation of a power series rather than a 2-term case) in a similar way that Gauss studied a hypergeometric series on the boundary of the disc of convergence.
The reason why we adopt an established technique (Gauss's Test) instead of new modern techniques for our research is that it is the only way to have appropriate values of the function within the radius of convergence (including its boundary), without violating a uniqueness theorem.  
\begin{thm} 
Let $y(x)= x^{\lambda }\left( d_0 +d_1 x+ d_2 x^2+ d_3 x^3+ \cdots \right)$, $d_i \ne 0$, be a power series (including a local Heun function) which consists of the 3-term recurrence relation. We assume that the polynomials in the numerator and denominator of $A_n$ and $B_n$ have the same degree:
\begin{equation}
d_{n+1}=A_n \;d_n +B_n \;d_{n-1} \hspace{1cm};n\geq 1 \label{eq:77}
\end{equation} 
\begin{equation}
\begin{cases} 
A_n =  \dfrac{C_t n^t + C_{t-1}n^{t-1} + \cdots + C_0}{c_t n^t + c_{t-1}n^{t-1} + \cdots + c_0}  = A  \dfrac{ n^t + \Omega _{t-1}n^{t-1} + \cdots + \Omega _0}{ n^t + \omega _{t-1}n^{t-1} + \cdots + \omega _0} \\
B_n = \dfrac{G_t n^t + G_{t-1}n^{t-1} + \cdots + G_0}{g_t n^t + g_{t-1}n^{t-1} + \cdots + g_0}  = B \dfrac{ n^t + \Theta _{t-1}n^{t-1} + \cdots + \Theta _0}{ n^t + \theta _{t-1}n^{t-1} + \cdots + \theta _0} \\
d_1= \dfrac{C_0}{c_0}  \;d_0 = A \dfrac{\Omega _0}{\omega _0}  \;d_0   
\end{cases}\nonumber
\end{equation}
where neither $C_t$ nor $c_t$ is zero and neither $G_t$ nor $g_t$ is zero with $\Omega _j= C_j/C_t$, $\omega _j = c_j/c_t $, $\Theta _j = G_j/G_t $ and $\theta _j= g_j/g_t$. Also, we denote $A=C_t/c_t$ and $B= G_t/g_t$.
 
The domain of absolute convergence of $y(x)$ is written by \cite{Choun2013} 
\begin{equation}
\mathcal{D} :=  \left\{  x \in \mathbb{C} \Bigg|  \left| \lim_{n\rightarrow \infty} A_n x \right| +\left| \lim_{n\rightarrow \infty} B_n x^2 \right|=  \left| A x \right| +\left| B x^2 \right| <1 \right\}
 \nonumber
\end{equation}
If $ \left| A x \right| +\left| B x^2 \right| =1$, then the series cannot converges. 





\end{thm} 
\begin{pf}
\textbf{1. The case of} $\mathbf{\Omega _{t-1}< \omega _{t-1}}$ \textbf{and} $\mathbf{ \Theta _{t-1} < \theta _{t-1}}$ 

We define  
\begin{equation} 
\left| \overline{A}_n \right| =  \left|  \dfrac{ n^t + \Omega _{t-1}n^{t-1} + \cdots + \Omega _0}{ n^t + \omega _{t-1}n^{t-1} + \cdots + \omega _0} \right|  
 \nonumber
\end{equation}
We also define $P(n)= n^t + \Omega _{t-1}n^{t-1} + \cdots + \Omega _0$ and $p(n)= n^t + \omega _{t-1}n^{t-1} + \cdots + \omega _0$.
 
For a large value of $n$, the leading coefficient determines the sign of the polynomial. After passing the rightmost root of a polynomial, a polynomial takes only positive values if  the leading coefficient is positive \cite{Bress}. It follows that once $n$ is larger than the largest root of either $P(n)$ or $p(n)$, then
\begin{equation} 
\left| \overline{A}_n \right| = \left|\frac{P(n)}{p(n)} \right| =  \frac{P(n)}{p(n)}
 \nonumber
\end{equation}  
If $\Omega _{t-1}< \omega _{t-1}$, then we can find a positive integer $h_1$ such that 
\begin{equation}
\omega _{t-1} - \Omega _{t-1}- h_1 <0 
\nonumber
\end{equation}
We observe that 
\begin{equation}
\frac{n}{n-h_1} \left| \overline{A}_n \right| = \dfrac{ n^{t+1} + \Omega _{t-1}n^{t} + \cdots }{ n^{t+1} + \left( \omega _{t-1}- h_1\right) n^{t} + \cdots } 
\nonumber
\end{equation} 
Since $ \omega _{t-1}- h_1 < \Omega _{t-1}$, the last fraction will eventually be larger than 1. We can say
\begin{equation}
 \left| \overline{A}_n \right| > 1-\frac{h_1}{n}
\nonumber
\end{equation}
Also, given positive error bound $\epsilon $ and a positive integer $N$,
$n\geq N$ implies that 
\begin{equation}
 \left| \overline{A}_n \right| > 1-\frac{h_1}{n} > 1-\epsilon 
\label{eq:5}
\end{equation}
We define 
\begin{equation} 
\left| \overline{B}_n \right| =  \left|  \dfrac{ n^t + \Theta _{t-1}n^{t-1} + \cdots + \Theta _0}{ n^t + \theta _{t-1}n^{t-1} + \cdots + \theta _0} \right|  
 \nonumber
\end{equation}
If  $\Theta  _{t-1}< \theta _{t-1}$, then we can find a positive integer $h_2$ such that 
\begin{equation}
\theta _{t-1} - \Theta _{t-1}- h_2 <0 
\nonumber
\end{equation}
We observe that 
\begin{equation}
\frac{n}{n-h_2} \left| \overline{B}_n \right| = \dfrac{ n^{t+1} + \Theta _{t-1}n^{t} + \cdots }{ n^{t+1} + \left( \theta _{t-1}- h_2\right) n^{t} + \cdots } 
\nonumber
\end{equation} 
Similarly, we can say
\begin{equation}
 \left| \overline{B}_n \right| > 1-\frac{h_2}{n}
\nonumber
\end{equation}
And, given positive error bound $\epsilon $ and a positive integer $N$,
$n\geq N$ implies that  
\begin{equation}
 \left| \overline{B}_n \right| > 1-\frac{h_2}{n} > 1-\epsilon 
\label{eq:6}
\end{equation}

We denote $\overline{i}= N+i$ and $\hat{i} = N+1+i$ where $i\in \{ 0,1,2,\cdots \}$.
For $n=N, N+1, N+2, \cdots$ in succession, take the modulus of the general term of $d_{n+1}$ in  (\ref{eq:77})
\begin{equation} 
 |d_{N+j}|= |\overline{c}_{j}| |d_N| + |\hat{c}_{j-1}| |B_{\overline{0}}|  |d_{N-1}|   
  \label{eq:7}
\end{equation}
where $j\in \{1,2,3,\cdots \}$, and  we define $|\overline{c}_{0}|= |\hat{c}_{0}|=1$.

In (\ref{eq:7}) is $|\overline{c}_{n}|$ is the sequence of the 3-term recurrence relation such as 
\begin{equation} 
 |\overline{c}_{n+1}|= \left|A_{\overline{n}} \right|  |\overline{c}_{n}| + \left|B_{\overline{n}} \right|  |\overline{c}_{n-1}|,\;\; n\geq 1   
  \label{eq:10}
\end{equation}
where $ |\overline{c}_{1}|= \left|A_{\overline{0}} \right|$.

Similarly, $|\hat{c}_{n}|$ is the sequence of the 3-term recurrence relation such as
\begin{equation} 
 |\hat{c}_{n+1}|= \left|A_{\hat{n}} \right|  |\hat{c}_{n}| + \left|B_{\hat{n}} \right|  |\hat{c}_{n-1}|,\;\; n\geq 1   
  \label{eq:11}
\end{equation}
where $ |\hat{c}_{1}|= \left|A_{\hat{0}} \right|$.

According to (\ref{eq:7}), then the series of absolute values, $1+|d_1||x|+|d_2||x|^2+ |d_3||x|^3 +\cdots$, is dominated by the convergent series
\begin{multline}
\sum_{n=0}^{N-1}|d_n||x|^n + |d_N||x|^N + \left( \overline{c}_{1}| |d_N| +   |B_{\overline{0}}|  |d_{N-1}|  \right) |x|^{N+1}
+ \left( \overline{c}_{2}| |d_N| + |\hat{c}_{1}| |B_{\overline{0}}|  |d_{N-1}| \right) |x|^{N+2} \\
+\left( \overline{c}_{3}| |d_N| + |\hat{c}_{2}| |B_{\overline{0}}|  |d_{N-1}|  \right) |x|^{N+3}
+\left( \overline{c}_{4}| |d_N| + |\hat{c}_{3}| |B_{\overline{0}}|  |d_{N-1}|  \right) |x|^{N+4}+\cdots \\
= \sum_{n=0}^{N-1}|d_n||x|^n +   |d_N|  |x|^N \sum_{i=0}^{\infty }\left| \overline{c}_{i} \right||x|^{i} + |B_{\overline{0}}|  |d_{N-1}|  |x|^{N+1} \sum_{i=0}^{\infty }\left| \hat{c}_{i}\right| |x|^{i}
\label{eq:12}
\end{multline}
By rearranging coefficients $\left| A_{\overline{n}}\right|$ and $\left|B_{\overline{n}}\right|$ in each sequence $\left|\overline{c}_{n}\right|$ in (\ref{eq:10}), 
\begin{equation}
\sum_{i=0}^{\infty }\left|\overline{c}_{i}\right| |x|^{i} = \left| \overline{y}_0(z)\right| + \left|\overline{y}_1(z)\right| \eta + \sum_{\tau =2}^{\infty } \left|\overline{y}_{\tau}(z)\right| \eta ^{\tau}  \label{eq:13}
\end{equation}
where
\begin{align}
\left|\overline{y}_0(z)\right| &= \sum_{i_0=0}^{\infty }  \prod _{i_1=0}^{i_0-1}\left| \overline{B}_{2i_1+1+N} \right| z^{i_0}\nonumber\\
\left|\overline{y}_1(z) \right| &=  \sum_{i_0=0}^{\infty } \left| \overline{A}_{2i_0+N}\right| \prod _{i_1=0}^{i_0-1}\left|\overline{B}_{2i_1+1+N} \right| \sum_{i_2=i_0}^{\infty } \prod _{i_3=i_0}^{i_2-1}\left|\overline{B}_{2i_3+2+N}\right| z^{i_2}\nonumber\\
\left|\overline{y}_{\tau}(z)\right|  &=   \sum_{i_0=0}^{\infty } \left|\overline{A}_{2i_0+N}\right| \prod _{i_1=0}^{i_0-1} \left| \overline{B}_{2i_1+1+N}\right|
 \prod _{k=1}^{\tau -1} \left( \sum_{i_{2k}= i_{2(k-1)}}^{\infty } \left|\overline{A}_{2i_{2k}+k+N}\right| \prod _{i_{2k+1}=i_{2(k-1)}}^{i_{2k}-1}\left|\overline{B}_{2i_{2k+1}+k+1+N}\right|\right)  \nonumber\\
  & \hspace{.3cm}\times \sum_{i_{2\tau } = i_{2(\tau -1)}}^{\infty }  \prod _{i_{2\tau +1}=i_{2(\tau -1)}}^{i_{2\tau }-1} \left|\overline{B}_{2i_{2\tau +1}+\tau +1+N}\right|  z^{i_{2\tau }}    \label{eq:14}
\end{align}
and
\begin{equation}
\begin{cases}
\eta = |A| |x| \cr
z= |B| |x|^2 \cr
0<\eta <1, \;\;\; 0<z<1 \cr
\end{cases}\nonumber 
\end{equation}
 The sequence $\left|\overline{c}_{n}\right|$ combines into combinations of $\left|A_{\overline{n}}\right|$ and $\left|B_{\overline{n}}\right|$ terms in (\ref{eq:10}): (\ref{eq:13}) is done by letting $\left|A_{\overline{n}}\right|$ in the sequence $\left|\overline{c}_{n}\right|$ is the leading term in a series $\sum_{i=0}^{\infty }\left|\overline{c}_{i}\right| |x|^{i}$; we observe the term of sequence $\left|\overline{c}_{n}\right|$ which includes zero term of $\left|A_{\overline{n}}\right|'s$ for a sub-power series $\left|\overline{y}_0(z)\right|$, one term of $\left|A_{\overline{n}}\right|'s$ for the sub-power series $\left|\overline{y}_1(z)\right|$, two terms of $\left|A_{\overline{n}}\right|'s$ for a $\left|\overline{y}_2(z)\right|$, three terms of $\left|A_{\overline{n}}\right|'s$ for a $\left|\overline{y}_3(z)\right|$, etc.

Similarly, by rearranging coefficients $\left| A_{\hat{n}}\right|$ and $\left|B_{\hat{n}}\right|$ in each sequence $\left|\hat{c}_{n}\right|$ in (\ref{eq:11}), 
\begin{equation}
\sum_{i=0}^{\infty }\left|\hat{c}_{i}\right| |x|^{i} = \left| \hat{y}_0(z)\right| + \left|\hat{y}_1(z)\right| \eta +\sum_{\tau =2}^{\infty } \left|\hat{y}_{\tau }(z)\right| \eta ^{\tau }  \label{eq:15}
\end{equation}
here, $\left| \hat{y}_k(z)\right|$ where $k\in \{0,1,2,\cdots\}$ is given by replacing an index $N$ with $N+1$ in (\ref{eq:14}).

A formula for the asymptotic of the Pochhammer symbol $(a)_n$ as $n\rightarrow \infty $
\begin{equation}
(a)_n \sim \frac{\sqrt{2 \pi}}{\Gamma (a)} e^{-n} n^{n+a-1/2} \left( 1+ \mathcal{O}(1/n)\right)
\nonumber
\end{equation}
and
\begin{equation}
\frac{(a)_n}{(b)_n} \sim \frac{\Gamma (b)}{\Gamma (a)}  n^{a-b} \left( 1+ \mathcal{O}(1/n)\right), \hspace{1cm} n\rightarrow \infty 
\nonumber
\end{equation} 
then, we can say
\begin{eqnarray}
\frac{\left( \frac{2+r+N-h_2}{2} + i_{2r}\right)_{i_{2(r+1)}}}{\left( \frac{2+r+N}{2} + i_{2r}\right)_{i_{2(r+1)}}} 
&>& \frac{\Gamma \left( \frac{2+r+N}{2} + i_{2r}\right)}{2\Gamma \left( \frac{2+r+N-h_2}{2} + i_{2r}\right)}  \frac{1}{ i_{2(r+1)}^{h_2/2}} 
\label{eq:16}
\end{eqnarray} 
where $i_{2(r+1)} \geq m$, here $N$ and $h_2$ are some positive integers, $r\in \{ 0,1,2,\cdots\}$ and  $N-h_2>0$.

Putting  (\ref{eq:5}) and  (\ref{eq:6}) at $\eta \left|\overline{y}_1(z) \right|$ in (\ref{eq:14}), we take the following inequality such as
\begin{eqnarray}
\eta \left|\overline{y}_1(z) \right|  
&>& (1-\epsilon ) \eta \sum_{i_0=0}^{\infty } (1-\epsilon )^{i_0} \sum_{i_2=i_0}^{\infty } \frac{\left( \frac{2+N-h_2}{2}+i_0\right)_{i_2 -i_0}}{\left(  \frac{2+N}{2}+i_0\right)_{i_2 -i_0}} z^{i_2}   \nonumber\\
&>& (1-\epsilon ) \eta \sum_{i_0=m}^{\infty } (1-\epsilon )^{i_0} \sum_{i_2=i_0}^{\infty } \frac{\left( \frac{2+N-h_2}{2}+i_0\right)_{i_2}}{\left(  \frac{2+N}{2}+i_0\right)_{i_2}} z^{i_2}
 \label{eq:17}
\end{eqnarray}
Put $r=0$ in (\ref{eq:16}) and take the new (\ref{eq:16}) into (\ref{eq:17})
\begin{eqnarray}
\eta \left|\overline{y}_1(z) \right|  
&>& \frac{(1-\epsilon )\eta}{2}   \sum_{i_0=m}^{\infty } \frac{\Gamma \left(\frac{2+N}{2}+i_0\right)}{\Gamma \left(\frac{2+N-h_2}{2}+i_0\right)}(1-\epsilon )^{i_0} \sum_{i_2=i_0}^{\infty } \frac{z^{i_2}}{ i_2^{h_2/2}}  \nonumber\\
&=& \frac{(1-\epsilon )\eta}{2}   \sum_{i_0=m}^{\infty } \frac{\Gamma \left(\frac{2+N}{2}+i_0\right)}{\Gamma \left(\frac{2+N-h_2}{2}+i_0\right)}(1-\epsilon )^{i_0} \sum_{i_2=m}^{\infty } \frac{z^{i_2}}{ i_2^{h_2/2}}  -\cdots \nonumber\\
&=& \frac{(1-\epsilon )^{m+1}\eta}{2}   \frac{\Gamma \left(\tfrac{2+N+2m}{2} \right)}{\Gamma \left(\tfrac{2+N+2m-h_2}{2}\right)} \; _2F_1{\scriptstyle \left( 1, \tfrac{2+N+2m}{2}, \tfrac{2+N+2m-h_2}{2};1-\epsilon \right)} \sum_{k=m}^{\infty } \frac{z^{k}}{ k^{h_2/2}}  -\cdots  \label{eq:18}
\end{eqnarray}
Putting  (\ref{eq:5}) and  (\ref{eq:6}) at $\eta^2 \left|\overline{y}_2(z) \right|$ in (\ref{eq:14}),  the  inequality is followed as
\begin{eqnarray}
\eta^2 \left|\overline{y}_2(z) \right|  
&>& \left( (1-\epsilon ) \eta \right)^2 \sum_{i_0=0}^{\infty } (1-\epsilon )^{i_0} \sum_{i_2=i_0}^{\infty } (1-\epsilon )^{i_2 - i_0} \sum_{i_4=i_2}^{\infty } \frac{\left( \frac{3+N-h_2}{2}+i_2\right)_{i_4 -i_2}}{\left(  \frac{3+N}{2}+i_2\right)_{i_4 -i_2}} z^{i_4}   \nonumber\\
&>& \left( (1-\epsilon ) \eta \right)^2 \sum_{i_0=m}^{\infty } (1-\epsilon )^{i_0} \sum_{i_2=i_0}^{\infty } (1-\epsilon )^{i_2  } \sum_{i_4=i_2}^{\infty } \frac{\left( \frac{3+N-h_2}{2}+i_2\right)_{i_4  }}{\left(  \frac{3+N}{2}+i_2\right)_{i_4  2}} z^{i_4} 
 \label{eq:19}
\end{eqnarray}
Put $r=1$ in (\ref{eq:16}) and take the new (\ref{eq:16}) into (\ref{eq:19})
\begin{eqnarray}
\eta^2 \left|\overline{y}_2(z) \right|  
&>& \frac{\left( (1-\epsilon ) \eta \right)^2}{2} \sum_{i_0=m}^{\infty } (1-\epsilon )^{i_0} \sum_{i_2=i_0}^{\infty }\frac{\Gamma \left(\frac{3+N}{2}+i_2\right)}{\Gamma \left(\frac{3+N-h_2}{2}+i_2\right)} (1-\epsilon )^{i_2  } \sum_{i_4=i_2}^{\infty } \frac{z^{i_4}}{i_4^{h_2/2}} \nonumber\\
&=& \frac{\left( (1-\epsilon ) \eta \right)^2}{2} \sum_{i_0=m}^{\infty } (1-\epsilon )^{i_0} \sum_{i_2=m}^{\infty }\frac{\Gamma \left(\frac{3+N}{2}+i_2\right)}{\Gamma \left(\frac{3+N-h_2}{2}+i_2\right)} (1-\epsilon )^{i_2  } \sum_{i_4=m}^{\infty } \frac{z^{i_4}}{i_4^{h_2/2}} -\cdots \nonumber\\
&=& \frac{\left( (1-\epsilon )^{m+1} \eta \right)^2}{2 \epsilon }  \frac{\Gamma \left(\tfrac{3+N+2m}{2} \right)}{\Gamma \left(\tfrac{3+N+2m-h_2}{2}\right)} \; _2F_1{\scriptstyle \left( 1, \tfrac{3+N+2m}{2}, \tfrac{3+N+2m-h_2}{2};1-\epsilon \right)} \sum_{k=m}^{\infty } \frac{z^{k}}{ k^{h_2/2}}  -\cdots \hspace{1cm}
 \label{eq:20}
\end{eqnarray}
Putting  (\ref{eq:5}) and  (\ref{eq:6}) at $\eta^3 \left|\overline{y}_3(z) \right|$ in (\ref{eq:14}),  the inequality ensues such as
\begin{eqnarray}
\eta^3 \left|\overline{y}_3(z) \right|  
&>& \left( (1-\epsilon ) \eta \right)^3 \sum_{i_0=0}^{\infty } (1-\epsilon )^{i_0} \sum_{i_2=i_0}^{\infty } (1-\epsilon )^{i_2 - i_0} \sum_{i_4=i_2}^{\infty } (1-\epsilon )^{i_4 - i_2} \sum_{i_6=i_4}^{\infty } \frac{\left( \frac{4+N-h_2}{2}+i_4\right)_{i_6 -i_4}}{\left(  \frac{4+N}{2}+i_4\right)_{i_6 -i_4}} z^{i_6}   \nonumber\\
&>& \left( (1-\epsilon ) \eta \right)^3 \sum_{i_0=m}^{\infty } (1-\epsilon )^{i_0} \sum_{i_2=i_0}^{\infty } (1-\epsilon )^{i_2  } \sum_{i_4=i_2}^{\infty } (1-\epsilon )^{i_4  } \sum_{i_6=i_4}^{\infty } \frac{\left( \frac{4+N-h_2}{2}+i_4\right)_{i_6 }}{\left(  \frac{4+N}{2}+i_4\right)_{i_6 }} z^{i_6} 
 \label{eq:21}
\end{eqnarray}
Put $r=2$ in (\ref{eq:16}) and take the new (\ref{eq:16}) into (\ref{eq:21})
\begin{eqnarray}
\eta^3 \left|\overline{y}_3(z) \right|   
&>& \frac{\left( (1-\epsilon ) \eta \right)^3}{2}  \sum_{i_0=m}^{\infty } (1-\epsilon )^{i_0} \sum_{i_2=i_0}^{\infty } (1-\epsilon )^{i_2} \sum_{i_4=i_2}^{\infty }\frac{\Gamma \left(\frac{4+N}{2}+i_4\right)}{\Gamma \left(\frac{4+N-h_2}{2}+i_4\right)} (1-\epsilon )^{i_4  } \sum_{i_6=i_4}^{\infty }  \frac{z^{i_6}}{i_6^{h_2/2}} \nonumber\\
&=& \frac{\left( (1-\epsilon ) \eta \right)^3}{2}  \sum_{i_0=m}^{\infty } (1-\epsilon )^{i_0} \sum_{i_2=m}^{\infty } (1-\epsilon )^{i_2} \sum_{i_4=m}^{\infty }\frac{\Gamma \left(\frac{4+N}{2}+i_4\right)}{\Gamma \left(\frac{4+N-h_2}{2}+i_4\right)} (1-\epsilon )^{i_4} \sum_{i_6=m}^{\infty }  \frac{z^{i_6}}{i_6^{h_2/2}} -\cdots \nonumber\\
&=& \frac{\left( (1-\epsilon )^{m+1} \eta \right)^3}{2 \epsilon^2 }  \frac{\Gamma \left(\tfrac{4+N+2m}{2} \right)}{\Gamma \left(\tfrac{4+N+2m-h_2}{2}\right)} \; _2F_1{\scriptstyle \left( 1, \tfrac{4+N+2m}{2}, \tfrac{4+N+2m-h_2}{2};1-\epsilon \right)} \sum_{k=m}^{\infty } \frac{z^{k}}{ k^{h_2/2}}  -\cdots
 \label{eq:22}
\end{eqnarray}
Similarly, putting  (\ref{eq:5}) and  (\ref{eq:6}) at $\eta^4 \left|\overline{y}_4(z) \right|$ in (\ref{eq:14}),  the  inequality is followed as
\begin{equation}
\eta^4 \left|\overline{y}_4(z) \right| > \frac{\left( (1-\epsilon )^{m+1} \eta \right)^4}{2 \epsilon^3 }  \frac{\Gamma \left(\tfrac{5+N+2m}{2} \right)}{\Gamma \left(\tfrac{5+N+2m-h_2}{2}\right)} \; _2F_1{\scriptstyle \left( 1, \tfrac{5+N+2m}{2}, \tfrac{5+N+2m-h_2}{2};1-\epsilon \right)} \sum_{k=m}^{\infty } \frac{z^{k}}{ k^{h_2/2}}  -\cdots
 \label{eq:23}
\end{equation}  
By mathematical induction, we repeat this process and  construct inequalities of every $\eta^{\tau } \left|\overline{y}_{\tau }(z) \right|$ terms where $\tau  \geq  5$.
Substitute (\ref{eq:18}), (\ref{eq:20}), (\ref{eq:22}), (\ref{eq:23}) and including  inequalities of all $\eta^{\tau } \left|\overline{y}_{\tau }(z) \right|$ terms where $\tau  \geq  5$ into (\ref{eq:13})
\begin{eqnarray}
\sum_{i=0}^{\infty }\left|\overline{c}_{i}\right| |x|^{i} &>&  \frac{\epsilon }{2} \sum_{j=1}^{\infty } \frac{\Gamma \left(\tfrac{1+N+2m+j}{2} \right)}{\Gamma \left(\tfrac{1+N+2m-h_2+j}{2}\right)} \left( \frac{ (1-\epsilon )^{m+1} \eta }{\epsilon  }\right)^j  \; _2F_1{\scriptstyle \left( 1, \tfrac{1+N+2m+j}{2}, \tfrac{1+N+2m-h_2+j}{2};1-\epsilon \right)} \sum_{k=m}^{\infty } \frac{z^{k}}{ k^{h_2/2}}  -\cdots \nonumber\\
&>&  \frac{\left( 1-K\right) \epsilon }{2}  \sum_{j=1}^{\infty } \frac{\Gamma \left(\tfrac{1+N+2m+j}{2} \right)}{\Gamma \left(\tfrac{1+N+2m-h_2+j}{2}\right)}\left( \frac{ (1-\epsilon )^{m+1} \eta }{\epsilon  }\right)^j  \; _2F_1{\scriptstyle \left( 1, \tfrac{1+N+2m+j}{2}, \tfrac{1+N+2m-h_2+j}{2};1-\epsilon \right)} \sum_{k=m}^{\infty } \frac{z^{k}}{ k^{h_2/2}}\nonumber\\
&>&   \frac{\left( 1-K\right) \epsilon }{2}  \sum_{j=1}^{\infty } \frac{\Gamma \left(\tfrac{1+N+2m+j}{2} \right)}{\Gamma \left(\tfrac{1+N+2m-h_2+j}{2}\right)}\left( \frac{ (1-\epsilon )^{m+1} \eta }{\epsilon  }\right)^j  \sum_{k=m}^{\infty } \frac{z^{k}}{ k^{h_2/2}} \nonumber\\
&=&  \frac{\left( 1-K\right) \epsilon }{2} \left\{ -\frac{\Gamma \left(\tfrac{1+N+2m }{2} \right)}{\Gamma \left(\tfrac{1+N+2m-h_2 }{2}\right)} +  \frac{\Gamma \left(\tfrac{1+N+2m }{2} \right)}{\Gamma \left(\tfrac{1+N+2m-h_2 }{2}\right)} \; _2F_1 \left( 1, \tfrac{1+N+2m }{2}, \tfrac{1+N+2m-h_2 }{2};\left( \frac{ (1-\epsilon )^{m+1} \eta }{\epsilon  }\right)^2 \right)  \right.  \nonumber\\
&&+ \left. \frac{ (1-\epsilon )^{m+1} \eta }{\epsilon  } \frac{\Gamma \left(\tfrac{2+N+2m }{2} \right)}{\Gamma \left(\tfrac{2+N+2m-h_2 }{2}\right)}  \; _2F_1  \left( 1, \tfrac{2+N+2m }{2}, \tfrac{2+N+2m-h_2 }{2};\left( \frac{ (1-\epsilon )^{m+1} \eta }{\epsilon  }\right)^2 \right)   \right\} \sum_{k=m}^{\infty } \frac{z^{k}}{ k^{h_2/2}}
\hspace{1cm}\label{eq:24}
\end{eqnarray}
where $0<K<1$ and we know $z+\eta =1$ on the boundary of the disc of convergence. $\sum_{k=m}^{\infty } \frac{z^{k}}{ k^{h_2/2}} <\infty $ because of $0<z=   |B|\left( \frac{  -|A| +\sqrt{|A|^2+ 4|B|}}{2|B|}\right)^2  <1$. We know an error bound $\epsilon \ll 1$  and $m$ is a fixed positive integer. We can  arbitrarily set the value of $\epsilon $, so it becomes $\frac{ \epsilon }{ (1-\epsilon )^{m+1} } \ll 1$.  $_2F_1$ functions are divergent at 
$\frac{\epsilon}{  (1-\epsilon )^{m+1}}< \eta = \frac{|A|\left( -|A| +\sqrt{|A|^2+ 4|B|}\right)}{2|B|} <1$ in (\ref{eq:24}). So, $\sum_{i=0}^{\infty }\left|\overline{c}_{i}\right| |x|^{i} >\infty $.  Similarly, an inequality of $\sum_{i=0}^{\infty }\left| \hat{c}_{i}\right| |x|^{i}$ is taken by replacing $N$ with $N+1$ in (\ref{eq:24}). So, $\sum_{i=0}^{\infty }\left| \hat{c}_{i}\right| |x|^{i} >\infty$ at $\frac{\epsilon}{  (1-\epsilon )^{m+1}}< \eta = \frac{|A|\left( -|A| +\sqrt{|A|^2+ 4|B|}\right)}{2|B|} <1$.
Therefore, (\ref{eq:12}) is divergent if $ \Omega _{t-1}< \omega _{t-1} $ and $ \Theta _{t-1} < \theta _{t-1} $. 

\textbf{2. The case of} $\mathbf{\Omega _{t-1}\geq \omega _{t-1}}$ \textbf{and} $\mathbf{ \Theta _{t-1} \geq \theta _{t-1}}$  

If $ \Omega _{t-1} > \omega _{t-1} $, then $P(n)-p(n)$  will eventually be  positive:
\begin{equation}
\left| \overline{A}_n \right| = 1+ \frac{P(n)-p(n)}{p(n)}>1 \nonumber
\end{equation}
for all $n$ larger than the rightmost root of $p(n) \left( P(n)-p(n) \right)$. Then, we can say
\begin{equation}
\left| \overline{A}_n \right| > 1-\frac{h_3}{n}> 1-\epsilon \nonumber
\end{equation}
for some positive integer $h_3$ with given positive error bound $\epsilon $ where $n\geq N$.

Likewise, if $ \Omega _{t-1} = \omega _{t-1} $ and $ \Omega _{t-2} > \omega _{t-2} $, then $P(n)-p(n)$  will be  positive:
\begin{equation}
\left| \overline{A}_n \right|  >1 > 1-\frac{h_3}{n}> 1-\epsilon \nonumber
\end{equation}
for large $n$ and $n\geq N$.

If $ \Omega _{t-1} = \omega _{t-1} $ and $ \Omega _{t-2} < \omega _{t-2} $, then we can find a positive integer $h_4$ such that 
\begin{equation}
\omega _{t-2} - \Omega _{t-2}- h_4 <0 
\nonumber
\end{equation}
We observe that 
\begin{equation}
\frac{n^2}{n^2-h_4} \left| \overline{A}_n \right| = \dfrac{ n^{t+2} + \Omega _{t-1}n^{t+1} + \Omega _{t-2}n^{t}+ \cdots }{ n^{t+2} + \omega _{t-1} n^{t+1}+ \left( \omega _{t-2}- h_4\right) n^{t} + \cdots } 
\nonumber
\end{equation} 
Since $ \omega _{t-2}- h_4 < \Omega _{t-2}$, the last fraction will eventually be larger than 1. We can say
\begin{equation}
 \left| \overline{A}_n \right| > 1-\frac{h_4}{n^2} > 1-\frac{h_3}{n}>1-\epsilon 
\nonumber
\end{equation} 
for some positive integer $h_3$ with given positive error bound $\epsilon $ where $n\geq N$.

Similarly, if $ \Omega _{t-1} = \omega _{t-1} $, $ \Omega _{t-2} = \omega _{t-2} $, $\cdots$, $ \Omega _{t-q} = \omega _{t-q} $ and $ \Omega _{t-q-1} > \omega _{t-q-1} $ where $1\in \{ 1,2,3, \cdots\}$,
\begin{equation}
 \left| \overline{A}_n \right| >  1-\frac{h_3}{n}>1-\epsilon 
\nonumber
\end{equation}
If $ \Omega _{t-1} = \omega _{t-1} $, $ \Omega _{t-2} = \omega _{t-2} $, $\cdots$  and $ \Omega _{t-q-1} = \omega _{t-q-1} $,
\begin{equation}
 \left| \overline{A}_n \right| = 1 >  1-\frac{h_3}{n}>1-\epsilon 
\nonumber
\end{equation} 
If $ \Omega _{t-1} = \omega _{t-1} $, $ \Omega _{t-2} = \omega _{t-2} $, $\cdots$, $ \Omega _{t-q} = \omega _{t-q} $ and $ \Omega _{t-q-1} < \omega _{t-q-1} $,
\begin{equation}
 \left| \overline{A}_n \right| > 1-\frac{h_q}{n^{q+1}} > 1-\frac{h_3}{n}>1-\epsilon 
\nonumber
\end{equation} 
where $\omega _{t-q-1}- \Omega _{t-q-1} -h_q <0$.

So, if $\Omega _{t-1}\geq \omega _{t-1}$,
\begin{equation}
 \left| \overline{A}_n \right| >  1-\frac{h_3}{n}>1-\epsilon 
\label{eq:25}
\end{equation}
Likewise, if $\Theta _{t-1} \geq \theta _{t-1}$,
\begin{equation}
 \left| \overline{B}_n \right| >  1-\frac{h_4}{n}>1-\epsilon 
\label{eq:26}
\end{equation}
for some positive integer $h_4$ with given positive error bound $\epsilon $ where $n\geq N$.

(\ref{eq:5}) and (\ref{eq:25}) are in the same form of inequality. And (\ref{eq:6}) and (\ref{eq:26}) are also in the same form of inequality. We say that (\ref{eq:16}) is also satisfied with $h_4$. So, inequality of $\sum_{i=0}^{\infty }\left|\overline{c}_{i}\right| |x|^{i} $ is same as (\ref{eq:24}) by replacing $h_2$ with $h_4$ where $N-h_4>0$. So, $\sum_{i=0}^{\infty }\left|\overline{c}_{i}\right| |x|^{i} >\infty $ Similarly, $\sum_{i=0}^{\infty }\left| \hat{c}_{i}\right| |x|^{i} >\infty$. Therefore,  (\ref{eq:12}) is divergent if $ \Omega _{t-1}\geq \omega _{t-1} $ and $ \Theta _{t-1} \geq \theta _{t-1} $.

\textbf{3. The case of} $\mathbf{\Omega _{t-1}\geq  \omega _{t-1}}$ \textbf{and} $\mathbf{ \Theta _{t-1} < \theta _{t-1}}$  

Inequality of $\left| \overline{A}_n \right|$ is (\ref{eq:25}) and inequality of $\left| \overline{B}_n \right|$ is (\ref{eq:6}). So, inequality of $\sum_{i=0}^{\infty }\left|\overline{c}_{i}\right| |x|^{i} $ is same as (\ref{eq:24}). Then, $\sum_{i=0}^{\infty }\left|\overline{c}_{i}\right| |x|^{i}, \sum_{i=0}^{\infty }\left| \hat{c}_{i}\right| |x|^{i} >\infty $. Therefore,  (\ref{eq:12}) is divergent if $ \Omega _{t-1}\geq  \omega _{t-1} $ and $ \Theta _{t-1} < \theta _{t-1} $.

\textbf{4. The case of} $\mathbf{\Omega _{t-1}< \omega _{t-1}}$ \textbf{and} $\mathbf{ \Theta _{t-1} \geq  \theta _{t-1}}$ 

Inequality of $\left| \overline{A}_n \right|$ is (\ref{eq:5}) and inequality of $\left| \overline{B}_n \right|$ is (\ref{eq:26}). So, inequality of $\sum_{i=0}^{\infty }\left|\overline{c}_{i}\right| |x|^{i} $ is same as (\ref{eq:24}) by replacing $h_2$ with $h_4$. Then,  $\sum_{i=0}^{\infty }\left|\overline{c}_{i}\right| |x|^{i}, \sum_{i=0}^{\infty }\left| \hat{c}_{i}\right| |x|^{i} >\infty $. Therefore,  (\ref{eq:12}) is divergent if $ \Omega _{t-1}< \omega _{t-1}$ and $ \Theta _{t-1} \geq  \theta _{t-1} $.

Finally, we conclude that the power series consisting of the 3-term recurrence relation, including a local Heun function, does not converge on the boundary of the disc of convergence.
\qed

\end{pf}


\end{document}